\magnification 1100
\def\UDF{1.1}

\def\BVA{2.1}
\def\BRA{2.2}
\def\BCV{2.3}
\def\TOP{2.4}
\def\BD{2.5}
\def\JEI{2.6}
\def\EV{2.7}

\def\UP{3.1}
\def\UMM{3.2}
\def\UPS{3.3}
\def\HUP{3.4}
\def\PUP{3.5}

\def\FIB{4.1}
\def\UM{4.2}
\def\HUM{4.3}
\def\PUM{4.4}

\font\tenmsb=msbm10
\font\sevenmsb=msbm10 at 7pt
\font\fivemsb=msbm10 at 5pt
\newfam\msbfam
\textfont\msbfam=\tenmsb
\scriptfont\msbfam=\sevenmsb
\scriptscriptfont\msbfam=\fivemsb
\def\Bbb#1{{\fam\msbfam\relax#1}}

\def\sqr#1#2{{\vcenter{\vbox{\hrule height.#2pt \hbox{\vrule width.#2pt
height #1pt \kern #1pt \vrule width.#2pt}\hrule height.#2pt}}}}
\def\qed{\sqr74}

\def\mapdown#1{\Big\downarrow \rlap{$\vcenter{\hbox{$\scriptstyle#1$}}$}}
\def\mapup#1{\Big\uparrow \rlap{$\vcenter{\hbox{$\scriptstyle#1$}}$}}
\def\lmapdown#1{\llap{$\vcenter{\hbox{$\scriptstyle#1$}}\Big\downarrow $}}
\def\mapright#1{\smash{\mathop{\longrightarrow}\limits^{#1}}}

 \def\la{\longrightarrow}

 \def\ni{\noindent}
 \def\cl{\centerline}
\def\df{\noindent {\it Definition.\  }}
 \def\rk{\noindent {\it Remark.\  }}
 
 \def\pf{\noindent {\it Proof.\  }}

\def\d{\delta }
 
 \def\l{\lambda}
 \def\s{\sigma}
\def\a{\alpha}
\def\b{\beta}
\def\g{\gamma}

\font\gothicf=eufm10
\font\sgothic=eufm7
\font\ssgothic=eufm5
\textfont5=\gothicf
\scriptfont5=\sgothic
\scriptscriptfont5=\ssgothic

\def\C{{\Bbb C}}

\def\P{{\Bbb P}}
\def\Q{{\Bbb Q}}
\def\B{{\cal B}}

\def\W{{\cal W}}

\def\X{{\cal X}}
\def\Z{{\Bbb Z}}

\def\Aut{{\rm Aut}}
\def\Pic{{\rm Pic}}

\def\Hom{{\rm Hom}}
\def\Isom{{\rm Isom}} 
\def\Im{{\rm Im}}

\def\dual{\omega _f}


\def\Agqs{{\tt A}(g,q,s)}
\def\DC{{\tt D}_{\C}(g)}
\def\Dg{{\tt D}(g)}

\def\Mgqs{{\tt M}(g,q,s)}

\def\Pgqs{{\tt P}(g,q,s)}
\def\Pghs{{\tt P}(g,h,s)}
\def\Mghs{{\tt M}(g,h,s)}
\def\Pso{{\tt P}^0(g,q,s)}

\def\QK{{\rm Q}_K(C)}
\def\QC{{\rm Q}_{\C}(C)}
\def\QX{{\rm Q}_K(X)}

\def\g{{\tt g}(g)}
\def\q{{\tt q}(g,q,s)}

\def\Hgdrs{{\tt H}(g,d,r,s)}
\def\Ngdrs{{\tt N}(g,d,r,s)}

\def\FBS{F_g(B,S)}
\def\FVT{F_g(V,T)}
\def\FFB{F_{g'}(B',S')}
\def\HM{\Hom _Z \bigl[ \B ,\Mg \times _Z Z\bigr]}
\def\HMf{{\cal H}om _Z \bigl[ \B ,\Mg \times _Z Z\bigr]}
\def\HJ{\Hom _Z \bigl[ (\B ,{\underline \s}) ,{\cal J}^r\bigr]}
\def\HJe{\Hom _Z ^e\bigl[ (\B ,{\underline \s}) ,{\cal J}^r\bigr]}

\def\MB{M[B,g]}
\def\MS{M\bigl[(B,S),g\bigr]}

\def\M{ M_Z[({\cal B},{\underline \s} ),g]}

\def\Ms{\overline M[B,g]}
\def\MBS{\overline M[(B,S),g]}

\def\Fp{F_{\psi}}

\def\Hilb{{\rm Hilb}^{p(x)}[\P ^r]}
\def\Mg{{\overline M_g}}
\def\Mh{M_{h(x)}}
\def\MgN{{\overline M^N_g}}
\def\CgN{{{\cal C}^N_g}}

\cl{\bf ON CERTAIN  UNIFORMITY PROPERTIES}
\cl{\bf OF
CURVES OVER FUNCTION FIELDS }

\

\ni LUCIA CAPORASO

\ni Department of Mathematics,
Massachusetts Istitute of Technology,
Cambridge MA 02139 USA

\ni
caporaso@math.MIT.edu
 
\

\

\ni
{\bf Contents}

\ni
1. Introduction

\ni
2. Boundedness results for stable families

\ni
3. Uniform versions of the theorem of Parshin-Arakelov

\ni
4. Uniformity of rational points over function fields

\

\ni
{\bf Abstract.} A uniform version of the Shafarevich Conjecture 
for function fields (Theorem of Parshin-Arakelov) 
is proved, together with a uniform version of the geometric Mordell conjecture
(Theorem of Manin). Such results are generalized to base varieties 
 of arbitrary
dimension (i.e. function fields of arbitrary transcendence degree).

\

\

\cl{{\bf 1. Introduction}}

\

\ni
{\it The problem.} 
Let $B$ be a non-singular, projective, complex curve and let $S$ be a finite
subset of points of $B$. The theorem below was conjectured by Shafarevich and
proved by Parshin (assuming $S=\emptyset$) and Arakelov (general case) in [P]
and [A].
\proclaim Theorem P-A.
Fix $g\geq 2$.
There exist only finitely many non-isotrivial families of smooth curves of genus g over $B-S$.

Some motivation for the Shafarevich conjecture was that, 
as Parshin shows in [P], 
the Mordell
conjecture for function fields
(Theorem M below, which had already been proved
by Manin [Ma]) can be obtained as a corollary of the above Theorem P-A.
Moreover, the same implication holds true for the number field analog of the
conjectures of Shafarevich and Mordell.

\proclaim Theorem M.
Let $K$ be a function field and let $X$ be a non-isotrivial curve of genus at least $2$ defined over $K$. Then $X$ has finitely many
$K$-rational points.

Denote by $X(K)$ the set of $K$-rational points of $X$. 
Let $V$ be a smooth variety over $\C$ and let $T$ be a closed subset of $V$.
Define
$\FVT$ as the set of  
equivalence classes of non-isotrivial families of smooth curves
of genus
$g\geq 2$ over $V-T$ (see below).

The problem studied in this paper is how
$|\FVT |$ and $|X(K)|$ vary as $V$, $T$, $X$ and $K$ vary. 
The case where $V=B$ is a curve and $K=\C (B)$ is its function field is treated
first, to obtain uniform versions of Theorem P-A and Theorem M. 
More precisely it will be proved that $|\FBS |$ and $|X(K)|$ are uniformly bounded by a function of $g$, the genus of $B$ and the
cardinality of $S$ (Theorems \UP \  and \UM ). 

The generalizations to higher dimensional bases are obtained
by specializing to suitable curves in $V$.
We  obtain a uniformity result where the  varieties $V$ are quantified by
their degree in some fixed projective space (\PUP \  and \PUM). From this a more
intrinsic uniformity statement is obtained for so-called  ``canonically polarized
varieties", that is, smooth varieties $V$ with ample canonical bundle
$K_V$. A natural invariant  quantifying canonically polarized varieties is the
Hilbert polynomial $h(x)$ of the canonical polarization (i.e. $h(n) = \chi
(K_V^{\otimes n})$ ). In such a case the cardinalities of $\FVT$ and $X(K)$ are
bounded by a function of $g$, the Hilbert polynomial $h(x)$, and the canonical
degree of $T$ (see Theorems \HUP \ and
\HUM ).

The proof of the Shafarevich conjecture is made of two parts, usually called
``boundedness" and ``rigidity". We proceed like in [P], giving essentially a
relative version of Parshin's proof of the boundedness part. 
To obtain our uniform version, we combine this  with techniques from
 the modern theory of stable curves and their moduli. 
This paper is about families of curves; for the analog of the Shafarevich conjecture
for families of higher dimensional varieties, see [F] and the very recent preprint
[BV].

We prove the uniform version of Manin's Theorem as a corollary of such a uniform
version of the Theorem of Parshin-Arakelov.
Another  perspective on this issue was taken by Miyaoka in [Mi], where he
indicated how an effective version of the Mordell conjecture could be obtained from
the Sakai inequality. Our approach here is completely different.

It is interesting to compare  our results with  analogous problems in the arithmetic
setting. The uniformity conjecture for rational points of curves over number fields
remains open, despite much work on the subject, in recent years.

\proclaim Uniformity Conjecture.
Fix $g\geq 2$. There exists a number $N(g) $ such that for any number field $F$ 
there exist at most finitely many curves of genus $g$, defined
over $F$, and having more than $N(g)$ $F$-rational points.

\ni
This is often called the ``strong" Uniformity conjecture, 
to distinguish it from its weaker variant, where the number field is fixed. At
the moment, they are both open 
(nor was anyone able to  give evidence for their being false). 
Interest in such  types of problems was revived after it was proved in [CHM]
that the Uniformity conjectures are a consequence of the so-called Lang
diophantine conjectures (also open) about the distribution of rational points
on varieties of general type.

\

\ni 
{\it Preliminaries.} 
We work over $\C$.
Throughout the paper $g$, $q$ and $s$ will be non-negative integers such that
$g\geq 2$; $B$ will be a smooth, projective complex curve of genus $q$ and $S\subset
B$ a subset of cardinality $s$ (unless otherwise specified).
 
By a ``stable'' (respectively ``semistable'') curve we mean a curve which is 
stable (respectively semistable) in the sense of Deligne and Mumford [DM].
If $C$ is a semistable curve, we denote by $C^{st}$ the stable model of $C$.
$M_g$ (respectively $\Mg$) denotes the coarse moduli space of smooth (resp. stable) curves of genus $g$.
The divisor of $\Mg$ parametrizing singular curves is denoted by $\Delta$.
A nodal curve is said to have genus $g$ if its arithmetic genus is $g$.

If $f:X\la V$ is a morphism, we denote by $X_v$ the fiber of  $X$ over $v\in V$; if $\phi :X \la Y$ is a $V$-map,
we denote by $\phi _v :X_v \la Y_v$ the restriction of $\phi $ to the fiber over $v\in V$.
By a ``family of curves" we mean a projective, flat morphism  $X\la V$  of integral varieties such that $X_v$ is a curve for every
$v\in V$, smooth for $v$ in some dense, open subset of $V$. A family of curves is called ``isotrivial" if all of its smooth fibers
are isomorphic.
 A ``family of stable curves of genus $g$" is a family of curves $X\la V$  all of whose fibers are stable curves of genus
$g$.

Given two families of  curves $f:X\la U$ and $h:Y\la U$ we say that they are
``equivalent" if there is an automorphism $\alpha $ of $U$ and an isomorphism
$\epsilon :X\la Y$ such that $h\circ \epsilon = \alpha \circ f$. 
This notion is really necessary only in the case $\dim U >1$, where $U$ would
typically be $V-T$ (see the introduction). 

A fibration is a proper morphism $f:X\la V$ such that $X$ and $V$ are smooth integral varieties and the fiber
$X_v$ is smooth for $v$ in some open dense subset of $V$. A semistable fibration is a fibration all
of whose fibers are semistable curves. Given a semistable fibration $f:X\la V$ its stable model is the family
$Y\la V$ of stable curves such that for every $v\in V$ $Y_v = X_v^{st}$. $Y$ is uniquely determined by the
mapping $X \la \P (f_*\dual ^{\otimes n})\la B$ for $n\geq 3$ (where $\dual $ is the relative dualizing sheaf).

Given a polynomial $p(x)$ we denote by $\Hilb$ the Hilbert scheme parametrizing closed subschemes of $\P ^r$
having Hilbert polynomial equal to $p(x)$.  Given a scheme $Z$ 
and a scheme $Y\la Z$ quasi-projective and flat over $Z$, we denote by ${\rm Hilb}_Z[Y] $ the Hilbert scheme
parametrizing closed, flat $Z$-subschemes of $Y$. Given  $X\la Z$ projective and flat over $Z$, the functor
${\cal H}om _Z[X,Y]$ of $Z$-morphisms is represented by a scheme $\Hom _Z[X,Y]$
which is described as an open subscheme of ${\rm Hilb}_Z[X\times _ZY] $.

Let $X$ be a curve of genus $g$, $B$ a curve of genus $q$ and $S\subset B$ a set of cardinality $s$; let 
${\rm f}(X,B,S)$ be a numerical function of $X$,$B$ and $S$. 
It will be convenient to use the following terminology.

\

\df \  
 $\rm f$ is ``uniformly bounded" if
there exists a constant ${\tt c}(g,q,s)$ such that for all $X$, $B$, $S$ as above we have ${\rm f}(X,B,S)\leq {\tt
c}(g,q,s)$. We say that a (finite) set is uniformly bounded if its cardinality is uniformly bounded.

\
 
\ni
An example of a uniformly bounded set which is relevant to this paper is the following.
Let $e(g)$ be some positive integer depending on $g$, let $C(B,S,e(g))$ be the set of finite coverings of $B$ having degree at
most equal to $e(g)$, and ramified only over $S$. Then $C(B,S,e(g))$  is a uniformly bounded set
(this would fail in positive characteristic).
The claim, of course, follows from  Riemann's theorem and from the fact that an
element in $C(B,S,e(g))$  is given by a group homomorphism from the
fundamental group of $B\setminus S$ to the symmetric group on $e(g)$ elements.

Let $K$ be a field containing $\C$, let $C$ be a smooth, irreducible curve defined over $K$.
Define
$$
\QK :=\{(\phi ,D): D {\rm{\  smooth \  curve\   of\   genus \  }} \geq 2 \  {\rm{ over}}\ 
 \  K, \  \phi :C\to D\  \rm{   finite \  map}\}.
$$
So that if $g(C)\leq 1$, then $\QK$ is empty. If instead $g(C)=2$ then $\QK
=\Aut _K(C)$.
 We recall a famous  theorem about $\QC$:

\proclaim Theorem D.
$\QC$ is a finite set. More precisely,
there exists a number  $\DC$ such that for every curve $C\in M_g$ we have
$
|\QC |\leq \DC .
$

\ni
The first  assertion is the classical theorem of De Franchis ([D]). The second is a modern refinement, 
which actually comes  with an 
explicit bound for $\DC$.  See [AP] and [K] for current records and related references.
We observe that a non-effective proof of the  existence of a uniform bound $\DC$
(which is all  we need here), given the theorem of De Franchis, is a  simple application
of the existence of Hilbert schemes. 
We need a uniform version  over function fields.
\proclaim Proposition \UDF.
Fix $g\geq 2$. There exists a number $\Dg$ such that for any function field $K$ and for any
smooth curve $X$ of genus $g$ defined over $K$ we have
$
|\QX |\leq \Dg .
$

\rk  We shall see that $\Dg = \DC$. 

\pf
Let $(\phi ^i,Y^i)\in \QX$ be distinct elements, for $i=1,...,N$. By contradiction, suppose that
$N > \DC$. 
We need the following 

\ni{\bf Claim.}
{\it Let $B$ be an integral scheme, let $X$, $Y^1$ and $Y^2$ be families of stable
curves of genus at least $2$ over $B$.
For $i=1,2$, let $\phi ^i:X\la Y^i$ be  finite $B$-morphisms.
Then there exists an open subset $U$ of $B$ such that $\forall b\in U$, $\phi _b^1\neq \phi
_b^2$.}

\ni
If $\phi ^1 =\phi ^2$  we take $U=\emptyset $, hence we shall assume that the two morphisms
are distinct. 
The Lemma is also clear if $Y^1$ and $Y^2$ are not fiberwise isomorphic.
Suppose then that $\forall b\in B$ we have $Y^1_b \cong Y_b^2$. This does not imply that 
the total spaces are isomorphic. Suppose first that they are, so that
$Y^1\cong Y^2 \cong Y$
then $\phi ^i $ determines a map $\Phi^i:B\la \Hom _B [X, Y]$ defined by $\Phi^i(b)=\phi^i_b$.
Clearly $\phi _b^1=\phi _b ^2$ if and only if $\Phi^1(b)=\Phi^2(b)$, since
$ \Phi ^1(B)\cap \Phi ^2(B)$ is closed we are done.

Now the general case. Given that $Y^1$ and $Y^2$ are fiberwise isomorphic,
we can construct a natural covering of $B$ where they become isomorphic.
This will be defined by the scheme over B
$
\Isom _B[Y^1,Y^2] \la B,
$
representing the functor of $B$-isomorphisms of the two families.
It is well known that 
$
\Isom _B[Y^1,Y^2] \la B
$
is a finite and unramified map (cf [DM]). Let $B'$ be a connected component dominating $B$ and let
$X', Z^1 ,Z^2 $ be the base changed families of $X$ $Y^1$, and $Y^2$ respectively.
Now $Z^1\cong Z^2 $ over $B'$.
Let
$\chi ^i:X'\la Z^i$ for $i=1,2$ be the pull-back of $\phi ^i$.
By the previous discussion, there exists an open subset $U'\subset B'$ where 
$\chi ^1 $ and $\chi ^2$ have distinct restrictions on the fibers. Since $B'$ is finite over
$B$, the image of $U'$ in $B$ is open, and we are done.

Let then $U_{i,j}$ be the open subset of $B$ such that
$\forall b\in U_{i,j }$ we have $\phi ^i_b\neq \phi _b^j$. For $i\neq j$ we have that
$U_{i,j}$ is not empty.
Let 
$U:=\bigcap _{1\leq i<j\leq N}U_{i,j},
$
of course $U$ is open and non empty.
Let $b\in U$, then by construction we have $N$ distinct elements $(\phi _b^i,Y_b^i)$ in 
$Q_{\C}(X_b)$. But $N>\DC$, while  by the refinement of the theorem of  De Franchis
 over the complex numbers we know that $|Q_{\C}(X_b)|\leq \DC$. We reached a contradiction.
$\qed$

\

\cl{\bf{2. Boundedness results for
stable families}}
 
\

The following set up will be fixed throughout the section.
Let $f:X\la B$ be a non-isotrivial  family
of stable curves of genus
$g$. The surface
$X$ needs not be smooth:  it will be normal with  isolated singularities  of type $xy=t^n$ at the
nodes of the fibers. We shall consider the relatively minimal desingularization $r:Y\la X$ of $X$,
with  its natural map
$h = f\circ r:Y\la B$. Then $h:Y\la B$ is a semistable fibration.
We shall denote by $S\subset B$ the set parametrizing singular fibers (of both $f$ and $h$), and we denote $s:=|S|$.

A large part of the material of this section is adapted  from  [P] and [S].
Recall that, by Arakelov's  Theorem,  $\dual $ is an ample line bundle  (since the family is not isotrivial). It  easily
follows that
$\omega _h$ is a big and nef line bundle on $Y$.

All the results of this section hold uniformly for any smooth complex curve $B$ of genus $q$ and  for any stable, non-isotrivial
family
$f:X\la B$ of curves of genus $g$ having at most $s$ singular fibers. We shall not repeat this in every statement.

\proclaim Lemma \BVA .
Let $g\geq 2, q, e$ be fixed integers.
There exists a number ${\tt n}(g,q,e)$ such that for any line bundle $E\in \Pic B$ with
$\deg E=e$ and for any $n\geq {\tt n}(g,q,e)$:

\ni
(1) $\dual ^{\otimes n} \otimes f^*E $ {\it is ample};

\ni
(2) $\omega _h ^{\otimes n}\otimes h^*E $ {\it is nef and big.}

\
 
\rk
The bound above does not depend on the number of singular fibers.

\pf
Clearly (2) follows from (1).
To prove (1) we use the Nakai Moishezon criterion (which holds for a normal surface).
Denote $L_n =\dual ^{\otimes n} \otimes f^*E $; it is easy to see  that if $e\geq 0$ then $L_n$ is ample for all $n\geq 1$,
hence we can assume that $e<0$.
We have
$$
L_n^2 = n^2\dual ^2+2ne(2g-2)
$$
and, since $\dual ^2 \geq 1$, if 
$
n >  -2e(2g-2)
$
we have that $L_n^2 >0$.

Let now $C$ be an integral curve in $X$, then $\dual \cdot C >0$ 
because $\dual $ is ample. If $C$ is contained in a fiber,
$
L_n \cdot C=\dual ^{\otimes n} \cdot C >0.
$

Suppose now  that $f_{|C} :C \la B$ is a finite map.
We need the fact that there exists ${\tt n}'(g,q,e)$ such that if $m \geq {\tt n}'(g,q,e)$, then $L_m$ is effective.

Before proving this, let us see how it implies the Lemma.
Let $L_m = {\cal O}(\sum r_i C_i)\otimes {\cal O}(V)$ such that $V$ is an effective divisor supported on the fibers of
$f$, $r_i >0$  and $C_i$ are irreducible curves  that cover $B$ via $f$.

If $C\neq C_i$, then $L_m \cdot C\geq 0$, as sum of intersections of effective, irreducible divisors with distinct supports;
therefore $L_{m+h} \cdot C>0$ for every $h\geq 1$.

If $C=C_i$, then 
$
\deg \{ f_{|C} :C\la B \} \leq m(2g-2),
$
hence 
$$
L_m \cdot C  \geq m(\dual \cdot C) + e m(2g-2)>e m(2g-2).
$$
Therefore if $h>  - em(2g-2) $ we have that $L_{m+h} \cdot C >0$.
It remains to prove that there is ${\tt n}'(g,q,e)$ such that $L_m$ is effective for $m \geq {\tt n}'(g,q,e)$.

For $n\geq 2$, $R^1f_*L_n =0$ hence
$$
\chi (L_n) = \chi (f_* L_n) = {n \choose 2}(12 \l - \d) + \l - (2n-1)(-e+q -1)(g-1) 
$$
where we used the formula $\deg f_* \dual ^n = {n \choose 2}(12 \l - \d) + \l$ for $n\geq 2$ proved in [Mu]
(see [Mu] for the definition of $\lambda $ and $\d$).
Then, if $q\neq 0$,
$$
\chi (L_n)\geq {n ^2 \over 2}(12 \l -\d) - n \bigl( {12 \l - \d \over 2} + 2(q-e+1)(g-1)\bigr)
$$
(recall that we are assuming $e<0$).
We conclude that if $n> 1+4(q-e+1)(g-1)$ and $q\geq 1$, we have
$\chi (L_n) >0$; we leave it to the reader to find the bound if $q=0$. 

It is easy to see that $h^2(X,L_n)=0$ (see the proof of \BCV (c) below)
and hence
$H^0(X,L_n)\neq 0$ for such a choice of $n$.
$\qed$

Let $\omega _X$ be the dualizing sheaf of $X$. We need a basic 
\proclaim Lemma \BRA.
Let $f:X\la B$ be a family of  stable curves. Let $A\in \Pic X$ be a relatively very ample line
bundle. Then $A$ is very ample on $X$ if, for any line bundle $M\in \Pic B$ such that $\deg M = -1,-2$, we have that
$A\otimes \omega _X^{-1}\otimes f^*M$ is nef and big.

\pf
It is well known that $A$ is very ample if and only if for every pair of (not necessarily distinct)
points $p$ and $q$ in $X$, we have
that
$
H^1(X,A\otimes {\cal I}_p)=H^1(X,A\otimes {\cal I}_p\otimes {\cal I}_q)=0.
$
Let $F_p$ and $F_q$ be the fibers of $f$ passing through $p$ and $q$ respectively.
We have  exact sequences:
$$
0\la A\otimes {\cal I}_{F_p}\la A\otimes {\cal I}_p\la  {\cal I}_p\otimes A_{|F_p}\la 0 $$

and
$$
0\la A\otimes {\cal I}_{F_p}{\cal I}_{F_q}\la A\otimes {\cal I}_p\otimes {\cal I}_q\la  {\cal I}_q{\cal I}_p\otimes A_{|F_p\cup F_q}\la 0.
$$
Since $A$ is relatively very ample, we have 
$
H^1(F_p,{\cal I}_p\otimes A_{|F_p})=H^1(F_p\cup F_q ,{\cal I}_q{\cal I}_p\otimes A_{|F_p\cup F_q})=0
$
for every $p$ and $q$ in $X$. Therefore, to show that $A$ is very ample,  it is enough to prove that
$$
H^1(X, A\otimes {\cal I}_{F_p})=H^1(X, A\otimes {\cal I}_{F_p}\otimes  {\cal I}_{F_q})=0,
$$
that is, it is enough to prove that for every line bundle $M$ on $B$ of degree equal to $-1$ or $-2$ 
we have 
$H^1(X, A\otimes f^*M)=0$.

By Serre's duality,
$H^1(X, A\otimes f^*M)\cong H^1(X, \omega _X\otimes (A\otimes f^*M)^{-1})^{-1}$.
We conclude by the Kawamata-Viehweg vanishing theorem, which can be applied to our surface $X$
(see [V] 2.35).
According to this theorem, if $L$ is nef and big, then $H^1(X,L^{-1})=0$;
we apply it to $L=  \omega _X^{-1}\otimes A\otimes  f^*M$, and conclude the proof.
$\qed$

\proclaim Corollary \BCV.
Fix $g\geq 2$ and $q$. There exists an integer ${\tt n}_0(g,q)$ such that for $n\geq {\tt
n}_0(g,q)$,

(a) $\dual ^n$ {\it is very ample on} $X$.

(b) $H^1(X, \dual ^n)=H^1(Y, \omega _h ^n)=0$.

(c) $H^2(X, \dual ^n)=H^2(Y, \omega _h ^n)=0$.

\

\rk (a) is in [S] (Lemme 9). We  include a proof as it takes two lines after what we already proved.
Observe also that from Matsusaka's big theorem  one derives that there exists a
constant $k_0$ depending on the Hilbert polinomial
$\chi (\dual ^n)$ such that $\dual ^n$ is very ample for $n\geq k_0$. This is a weaker version of $(a)$ below, as $k_0$ depends on
$g,q$ and $s$.

\ni
\pf By Lemma \BRA ,
$\dual ^n$ is very ample if $\dual ^{n-1} \otimes f^*(M\otimes K_B ^{-1})$ is nef and
big for every $M$ of degree $-1$ or $-2$ on $B$. Then we apply Lemma \BVA \  and
obtain (a).

(b) follows by Lemma \BVA  and Kawamata-Viehweg vanishing theorem (see above) applied to
the line bundle
$\dual ^{n-1}\otimes f^*K_B^{-1}$. Similarly for $\omega _h$. 

It is enough to prove (c) for $\omega _h$. By Serre duality

$$
H^2(X, \omega _h ^n)=H^0(X,\omega _h ^{\otimes (1-n)}\otimes h^*K_B)
$$
now $\omega _h ^{\otimes (1-n)}\otimes h^*K_B\cdot X_b<0$ for $n\geq 2$. Since $X_b ^2 =0$ we obtain that $\omega _h ^{\otimes
(1-n)}\otimes h^*K_B$ cannot possibly be effective.
$\qed$
 \

\proclaim Lemma \TOP. 
$\omega _h ^2 $ and $\chi ({\cal O}_Y )$ are uniformly bounded.

\pf The first  part of the statement is  in [S], Th\'eor\`eme 3, where more precisely
it is proven that $\omega _h ^2 \leq 8g(g-1)(q-1+{s\over 2})$. The uniform boundedness of 
$\chi ({\cal O}_Y )$ follows from [P] Proposition 1.
$\qed$

\

\proclaim Corollary \BD .
Fix $n$ so that \BCV \    holds. Then there exists a number $r_0={\tt r}_0(g,q,s)$ such that
$h^0(X, \dual ^n)\leq r_0$.

\pf It is sufficient to show the statement for $\omega _h ^n$.
We have by Riemann-Roch and by \BCV
$$
h^0(Y, \omega _h ^n) = \chi (\omega _h ^n) = \chi ({\cal O}_Y) + {n^2 \omega _h  ^2 - n (\omega _h \cdot K_Y) \over 2}
$$
and 
$(\omega _h \cdot K_Y) = \omega _h ^2 + \omega _h \cdot h^*K_B= \omega _h ^2 + (2g-2)(2q-2)  >(2g-2)(2q-2) .$
We conclude by the previous result. $\qed$

\

\ni
The above analysis can be now applied to perform a useful construction.
 Fix $n\geq
{\tt n}_0(g,q)$ so that \BCV \  holds, and let $d=n(2g-2)$.
Given $f:X\la B$, the line bundle $\dual ^n$ determines an embedding 
$$
\phi : X\hookrightarrow \P ^r =\P ( H^0(X, \dual ^n)^*)
$$
where $r\leq {\tt r}_0(g,q,s)$  as in Corollary \BD \  above. For every $b\in B$ the fiber $X_b$ is mapped isomorphically
by $\phi$
to a  curve spanning a linear subspace $\P ^k \subset \P ^r$, where $k=d-g$. In fact $H^1(X, \dual ^n\otimes {\cal I}_{X_b}) =0 $,
because we have chosen $n$ so that 
$\dual ^n\otimes {\cal I}_{X_b}$ is ample (see the proof of \BCV ).

Let $p(x)=dx - g +1$, fix $\P ^r$ as above and 
consider the Hilbert scheme $\Hilb$ with its universal curve ${\cal Y} \subset \Hilb \times \P ^r$; we now introduce 
the following subscheme of  $\Hilb$
$$
J^r := \{h\in \Hilb : Y_h \  {\rm{is\  stable\  and \  spans \  a}} \  \P ^k,\  \omega_{Y_h}^{\otimes n}
\cong {\cal O}_{\P ^k}(1)\otimes {\cal O}_{Y_h}\}.
$$
where  $Y_h$ denotes the fiber over $h\in \Hilb$ of the universal curve.
We denote by $\Delta ^r\subset J^r$ the locus of singular curves.

\proclaim Lemma \JEI .
$J^r$ is a smooth, irreducible, quasi-projective variety of dimension  $r(k+1)+k+3g-3$. $\Delta ^r$ is a reduced divisor with
normal crossings singularities.

\pf
If $k=r$ we denote  $J^r=H_g$; this is the subset of the Hilbert scheme parametrizing stable curves embedded by the
$n$-canonical linear series; since $r=d-g$ these curves are not contained in any hyperplane. The geometry of $H_g$ is well
understood, in fact the geometric invariant theory quotient of $H_g$ by the natural action of $SL(r+1)$ is identified with $\Mg$.
Results of [G] give our statement for $r=k$.
 If $k<r$ there is a natural surjective morphism 
$$
\l :J^r\la {\Bbb G} (k,r)
$$
where ${\Bbb G} (k,r)$ is the Grassmannian of $k$-dimensional linear subspaces in $\P ^r$.
The map $\l$ maps the  point $h$ to the linear span of the fiber $Y_h$ in $\P ^r$. Its fibers are all
isomorphic to $H_g$. To conclude just notice  that ${\Bbb G} (k,r)$ is a smooth and irreducible projective 
variety of dimension
$(k+1)(r-k)$.
$\qed$

\

\ni
The embedding $\phi :X\hookrightarrow \P ^r$ defined above defines a morphism $\Phi :B\la J^r$.
Fix a projective model $B\times J^r\subset \P ^s$ and let $\Gamma $ be the graph of $\Phi$,
so that $\Gamma \subset B\times J^r\subset \P ^s$.

\proclaim Lemma \EV .
There exists ${\tt e}_0(g,q,s)$ such that $\deg \Gamma \leq  {\tt e}_0(g,q,s)$.

\ni
This is essentially Lemma 5 in [P]. To adjust Parshin's proof to our situation one needs to consider heights with respect to the
polarization given by $\dual ^n$ rather than $K_X^6$.
Effective versions of such a result are given in [EV], for example,
or in the more recent [T] where other relevant references are given.

\
 
\cl{{\bf 3. Uniform versions of the theorem of Parshin-Arakelov}}

\

\ni
{\it One-dimensional base.} 
Fix $B$ a  curve and $S\subset B$ as usual. 

\

\df Let  $\FBS$ to be
 the set of equivalence classes of non-isotrivial fibrations  $f:X\to B$,
with $X$  a smooth, relatively minimal surface and, for all
$b\not\in B$, $X_b$ is a smooth curve of genus $g$.

\

\ni
By the existence and unicity of relatively
minimal models over $B$ for fibrations of genus $g$ over an open subset of $B$, this definition is consistent with the one
given in the introduction for $\FVT$.
The goal of this section is the following result

\proclaim Theorem \UP.
Fix $g\geq 2$, q and s.
There exists a number $\Pgqs$ such that for any  curve $B$ of genus q and for any
 subset $S\subset B$ of cardinality $s$, we have
$|\FBS |\leq \Pgqs$.

\ni
There are   few cases where the theorem is already known
(see [B]): if $q=0$ and $s\leq 2$ then $F_g(\P ^1, S) = \emptyset$ for every $g$;
if $q =1$ and $S=\emptyset$ again for every $g$, $F_g(B,\emptyset )=\emptyset$.
This fact  allows us to ignore the equivalence relation
for  families over $B$ (all automorphism groups being  uniformly bounded).

The theorem is proved in three steps. First we show that  for
every pair $(B,S)$ as above, the set of non-constant morphisms $\psi :B\la
\Mg$ such that $\psi ^{-1}(\Delta )\subset S$ and such that $\psi$ is the
moduli map (see below) for some family of stable curves over $B$, is uniformly
bounded. Then we show that for any fixed $\psi$ the set of families having
$\psi$ as moduli map  is uniformly bounded. Finally, we reduce the general
case to the stable case.

The first part  is a consequence of the construction of a moduli space for moduli
maps of curve of genus $q$ to $\Mg$.
Let $N$ be a coarse moduli space (for example, $N=M_g$ or $N=\Mg$)

\

\ni
\df A morphism $\psi :B\la N$ is a ``moduli map" if $\psi$ is non-constant and if there exists a family
$f:X\la B$ of objects parametrized by $N$ such that $\psi :b\mapsto [X_b]\  $ is
the  moduli map of $f$.

\

\ni
Clearly, if $N$ is a fine moduli space,  every non-constant map is a moduli
map. We denote 
$
\MB
$
(respectively $\Ms$) the set of moduli maps from $B$ to $N=M_g$ (respectively $N=\Mg$).
Let $S$ be a finite subset of $B$, denote by $\MBS$ the set of all moduli maps
$\psi $ of $B$ to $\Mg$ such that 
$\psi ^{-1}(\Delta )\subset S$. That is, $\MBS$ is the set of moduli maps with  degeneracy locus contained in $S$; it is a
finite set, by the theorem of Parshin-Arakelov.

\proclaim Proposition \UMM.
Fix $g\geq 2$, $q$, and $s$. There exists a number $\Pso$ such that for any  curve $B$ of genus
$q$ and any finite subset $S\subset B$ of cardinality $s$ we have
$|\MS |\leq \Pso .$
 
\pf
 In the relative setting, let $h:\B \la Z$ be a projective morphism of integral
 varieties
(we can assume it to be smooth).
We shall say that a $Z-$morphism $\Psi :\B  \la N\times Z$ is a moduli map over $Z$ if there exists a family $\X \la
\B $ such that $\forall z\in Z$ the restriction $\psi _z :B_z \la N\times \{ z\}$
is the moduli map for the family 
$X_z \la B_z$.

We are mainly interested in the following variant.
Given $h:\B \la Z$ as above,  consider the special case in which  $h$ is a family of
smooth curves, assume in addition that $h$ has $s$ non-intersecting sections
$\s _i :Z\la \B$ with $h\circ \s _i = {\rm id}_Z$ for every $i=1,....,s$.
Denote  $\underline \s = \{ \s _1,....,\s _s\}$, $\  {\cal S} = \cup \s _i (Z) \subset \B$, for $z \in Z$ denote by $S_z$ the
finite subset of $B_z$ given by $\{ \s _1(z),....,\s _s(z)\}$. Consider the $Z$-set $\M \la Z$ whose fiber over $z\in Z$ is 
the set $M[(B_z,S_z),g]$.

As a set, the above space of moduli maps with fiberwise fixed degeneracies is naturally a subset 
$$
\M \subset \HM .
$$
We are going to show that $\M$ is a finite union of irreducible, quasi-projective varieties; to do that we construct an
auxiliary parameter scheme
$A$ which has a natural morphism to $\HM$, whose image we shall identify with $\M$. 

Let $J^r$ be as defined at the end of the previous section  and
let ${\cal J}^r:=J^r\times Z$.
Define
$\HJ$ to be the subscheme  of $\Hom  _Z[\B ,{\cal J}^r]$ parametrizing morphisms $\Phi$ such that
$\Phi ^{-1} (\Delta ^r \times Z) \subset  {\cal S}$ and such that $\forall z \in Z$ the restriction
$\phi _z : B_z \la J ^r \times \{ z\}$ is not constant.
There is an inclusion of schemes
$$
\HJ \hookrightarrow {\rm Hilb } _Z\bigr[ \B \times _Z {\cal J}^r\bigl]
$$
associating to a morphism its graph.
Since $\B$ is projective over $Z$,
for a fixed integer  $e$  we can consider the quasi-projective scheme  over $Z$ parametrizing
morphisms of degree $e$, that is we consider 
$$
\HJe \hookrightarrow {\rm Hilb }^{ex-q+1} _Z\bigr[ \B \times _Z {\cal J}^r\bigl] .
$$
Let
$$
A:=\bigcup _{r\leq r_0}\bigcup _{e\leq e_0} \HJe
$$
where $e_0={\tt e}_0(g,q,r)$ and $r_0={\tt r}_0(g,q,r)$ are as in \EV \  and \BD \  respectively. $A$ is thus 
a scheme of finite type
over
$Z$. We now show that $A$ carries a natural family of moduli maps to $\Mg$, with fixed degeneracies.
This follows by the functorial properties of the various moduli spaces involved;
we have a diagram of objects an morphisms over $Z$:
\

$$
\matrix{&&A&\buildrel i \over\hookrightarrow&{\rm Hilb }_Z\bigr[ \B \times
_Z {\cal J} \bigl] &&\cr
&&\mapup {} & &\mapup {}&&\cr
&&{\cal G}&\mapright{j}&{\cal C} &\hookrightarrow & \B \times {\cal J} \times {\rm Hilb }_Z\bigr[
\B
\times _Z {\cal J} \bigl] \cr
& &\mapup {} & &\mapup {}& & \mapup{} &\cr
A\times \B &\longleftarrow &{\cal W} &\mapright{} &{\cal D} &\hookrightarrow & \B \times {\cal Y}
\times {\rm Hilb }_Z\bigr[ \B \times _Z {\cal J} \bigl] }
$$
\

\

\ni
where $\cal J$ is the disjoint union of ${\cal J}^r$, for $r\leq r_0$ ,  $\  {\cal Y} \la  {\cal J}$
is the universal family, and the Hilbert scheme is defined correspondingly. We denoted
${\cal G} = i^* {\cal C}$ and $\W =j^* {\cal D}$.
The above diagram yields 
$$
\matrix{\W &\mapright{}&A\times _Z\B&&\cr
\  \  \lmapdown {}& &\mapdown {}&&\cr
A&\mapright {} &Z &\longleftarrow &\B}
$$
so that for every $a\in A$, if $z$ is the image of $a$ in $Z$ and $B_z$ is the fiber of $\B$ over
$z$, the restriction
$W_a \la \{a\} \times B_z$ is a family of stable curves of genus $g$ with singular fibers over 
$(a, \s _i(z) )$ for $i=1,....,s$.
Therefore the above diagram can be viewed as an element of
$\HMf (A)$ and as such it determines a unique morphism
$$
\beta : A \la \HM .
$$
Now $\Im \b$ is  a finite union of quasi-projective schemes over $Z$, and, by construction, it is a subset of the set
of moduli maps with fixed degeneracies, that is
$$
\Im \b \subset \M .
$$
To show that  the above inclusion is in fact an
equality of sets we  apply the results of the previous section. Given a moduli map $\psi : B\la \Mg$, let $f:X\la B$ be the
corresponding family of stable curves. By the construction at the end of section 2 we can associate to such data a morphism
from $B$ to $J^r$ (canonically up to an action of the $PGL(r+1)$); this determines a unique
 point in $A\subset \HJ$ whose image via $\beta$ is the given $[\psi ]$.

To conclude the proof of the proposition, we shall
 apply the above construction to a suitable $h:\B \la Z$ with $s$ sections.
Suitable means that $Z$ must be some fine moduli space for curves of genus $q$ with $s$
marked points. There are various choices for this, for example,  $Z=M_{q,s}^N$, the moduli
space of smooth curves of genus $q$, with $s$ marked points, and level $N\geq 3$-structure.
For every pair $(B,S)$ as above there exists a (not necessarily unique) point $z\in
Z$ corresponding to such a pair (that is, $B_z=B$ and $\{ \s_1 (z),....,\s_n(z)\}
=S$). Now $\M$ is a finite union of irreducible quasi-projective varieties  over $Z$
and its fiber over the point $z$ is
$\MS$ (for every $z\in Z$).
By  the theorem of Parshin-Arakelov, $\M$ has finite fibers over $Z$; therefore there exists 
an upper bound   $\Pso$ on the cardinality of the fibers. $\qed$

\

\ni
This concludes the first part of the program to prove Theorem \UP . The discussion that follows is
needed in order to  handle curves with automorphisms.
We start by a general
fact about moduli maps over a general base scheme, which will be useful again later.
 Let  $V$ be a fixed quasi-projective variety of any dimension,
assume that
$V$ is irreducible and smooth, let $\psi :V\la \Mg$ be a morphism such that $\Im \psi \cap M_g \neq \emptyset$.
 
\

\df Let $\Fp$ be the set   of all equivalence classes of stable familes
$f:X\la V$ such that for every $v\in V$, $X_v$ is stable and  $\psi
(v) = [X_v]$.

\proclaim Lemma \UPS . For every $V$ and $\psi$ as above, 
$|\Fp |\leq \Dg$.

\pf 
Let $f^i :X^i \la V$, $i=1,....,n$,  be distinct elements of $\Fp$, and assume, by contradiction,
that $n>\Dg\  $  (see \UDF ). The families $f^i$ are fiberwise isomorphic. We 
construct a covering
$V'\la V$ over which they are all isomorphic to a certain family $X\la V'$. Let 
$V_i :=\Isom _V(X^i,X^n)$ so
that
$V_i$ has a natural, finite morphism to $V$ (Theorem 1.11 of [DM]). 
Replacing $V_i$ with a connected component dominating $V$ will not alter the argument.
Let $X^i_1 := X^i \times _V V_1$ for
$i=1,....,n$. Then
$X^1_1
\cong X^n_1$ over
$V_1$. Moreover, for each $i$ there is a natural commutative diagram
$$
\matrix{X^i_1&\mapright{}&X^i\cr
\  \  \lmapdown {f^i_1}& &\mapdown {f^i}\cr
V_1&\mapright {} &V}
$$
and over $V_1$ we have $n-1$ different familes $f^i_1$. Iterate this construction to get $V'\la V$.
We have a commutative diagram for every $i=1,...,n$ 
$$
\matrix{X&\mapright {\phi ^i} & Y^i &\mapright{}&X^i\cr
\mapdown {}& &\mapdown {}& &\mapdown {}\cr
V'&=&V' &\mapright {} &V}
$$
where $Y^i = X^i \times _V V'$ and $\phi ^i$ is an isomorphism.
By assumption,   $\phi ^i \neq \phi ^j$, we thus obtain $n>\Dg$ distinct elements in
$\QX$ ($K$ being the function field of $V'$), which is not possible, by \UDF . $\qed$

\

\ni
{\it Proof of Theorem \UP.}
Denote by $\  \FBS ^s$  the subset of $\FBS$ consisting of semistable fibrations. 
We   use
curves with level  structure to show that there is
a map
$$
\lambda :\FBS \la \cup F_g(B',S')^s
$$
where the union is for $(B',S')$ varying in a uniformly bounded set,
and the fibers are uniformly bounded sets.

 Fix $N\geq 3$ and let $\MgN$ be the moduli space of stable
curves of genus $g$ with level $N$ structure. This is a fine moduli space, endowed with a 
universal family $\CgN\la \MgN$ and with a finite morphism to $\Mg$. Let $[f:X\la B]\in \FBS$,
we have a commutative diagram
$$
\matrix{&&Z&\mapright{}&\CgN\cr
&&\  \  \lmapdown {}& &\mapdown h\cr
X' & \mapright{f'}&B'&\mapright {\psi '} &\MgN\cr
\mapdown{}& &\mapdown{\rho }& &\mapdown{} \cr
X&\mapright {f}& B & \mapright {\psi} & \Mg}
$$
where $f'$ is the rigidification of $f$ by level $N$ structures,  so that $\rho $ is a finite
morphism which is \'etale with group $Sp(g,\Z /n\Z )$ away from $S$.
Therefore 
$\rho$, which depends on $f$, varies in a uniformly bounded set as $f$ varies. Let $S' = \rho
^{-1}(S)$.

Away from singular fibers, $\psi '$ is the 
moduli map of $f'$, so that the pull back $Z$ of $\CgN$ to $B'$ is naturally birational to $X'$ over
$B'$ (in fact $Z$ is isomorphic to $X'$ away from $S'$). 
The family $Z\la B'$ is a stable reduction of $f$.

Now we define $\l$ by sending  the
fibration
$f:X\la B$ to the relatively minimal semistable fibration
obtained by resolving the singularities of the stable family $Z\la B'$. 
The previous discussion shows that the union defining the range of $\l$ is over $(B',S')$ varying in  a
uniformly bunded set; in addition,
the fibers of
$\l$ are  uniformly bounded, by  \UDF.

It remains  to show that $\FBS ^s$ is uniformly bounded.
There is a surjection

$$
\matrix {\FBS  ^s& \la &\MS \cr
f & \mapsto & \psi _f}
$$
where $\psi _f(b) = [X_b^{st}]$,
that is to say, $\psi _f$ is the moduli map of the stable model of $X\la B$.
 By  \UMM \    $\MS$ is
uniformly bounded. Since the smooth, relatively  minimal model is unique, the fiber of the above map over $\psi _f =\psi$ is
naturally in bijective correspondence  with  $\Fp$.
But $\Fp$ is uniformly bounded, by  \UPS  \    hence the Theorem is proved. $\qed$

\

\ni
{\it Higher dimensional base.}
We  apply the results of Theorem \UP \  
to families of  curves parametrized by  smooth,
projective, integral varieties $V$ of
arbitrary dimension.  
The first issue is how to quantify the base spaces, so as to obtain analog uniformity
statements. We proceed as follows: fix a polynomial $h\in \Q [x]$  such that $h(\Z
)\subset \Z$. Let $\Mh$ be the moduli space of smooth projective varieties $V$
having ample canonical line bundle $K_V$ and such that 
$
\chi (K_V^{\otimes n})=h(n)
$
It is well known (see [V]) that $\Mh $ is a quasi-projective scheme, called the moduli space of ``canonically polarized"
varieties. For example, if $\deg h =1$, in order for $\Mh$ not to be empty one needs that there exist an integer $q\geq 2$
such that
$h(x)=(2q-2)x -q+1$ and $\Mh$ is equal to the moduli space of curves of genus $q$.
We have the following generalization of Theorem \UP :

\proclaim Theorem \HUP.
Fix non-negative integers $g\geq 2$ and $s$ and fix a polynomial $h$ as above. There exists a number $\Pghs$ such that for any
$V\in \Mh$, for any closed $T\subset V$ such that $K_V^{\dim V -1}\cdot T\leq s$ 
there exist at most $\Pghs$ equivalence classes of
non-isotrivial families of smooth curves of genus $g$ over  $V-T$.

\pf
The proof follows from a more general uniformity statement for polarized varieties (\PUP \  below). First, we use Matsusaka's
big Theorem and fix an integer $\nu = \nu (h)$ depending on the Hilbert polynomial $h$ such that for every $n\geq \nu$ and for
every $V\in \Mh$, the line bundle $K_V^{\otimes n}$  is very ample. By the vanishing Theorem of Kodaira, it follows that
$h^0(V, K_V^{\otimes n}) = h(n)$ for $n>\nu$ and for every $V\in \Mh$.
Fix then  $n>\nu$  and let $r =h(n) -1 $;  we  fix a projective space $\P ^r$ such that for every $V\in \Mh$ 
the line bundle $K_V^{\otimes n}$ determines a projective
model of $V$ in $\P ^r$ having degree $d$ determined by $h$ (in fact, $d$ equals the leading 
coefficient of $h$, multiplied by
$\dim V !$). The Theorem follows from the following

\proclaim Lemma \PUP .
Fix $g$, $d$, $r$ and $s$ integers; there exists a number $\Hgdrs$ satisfying the following property.
 For any smooth projective variety
$V\subset \P ^r$ of degree $d$ 
and for any subvariety
$T\subset V$ of degree $s$, there exist at most $\Hgdrs$ 
equivalence classes of non-isotrivial families of smooth curves of genus $g$ over
$V-T$.

\pf
There exists a
number $q={\tt q}(d,r)$ such that $V$ is covered by smooth curves $B$ of genus  at most  $ q$ passing through any of its points.
In fact such 
$B$ are just one-dimensional linear sections of $V$, that is, they are obtained by intersecting $V$ with the correct number of
generic
hyperplanes in $\P ^r$. By a theorem of Castelnuovo (see [ACGH] III.2), the genus of a smooth curve of degree
$d$ in
$\P ^r$ is bounded above by a function ${\tt q}(d,r)$ of $d$ and $r$. In addition, any such $B$ intersects $T$ in at most
$s$ points. Let 
$$
P= \max _{q'\leq q, s'\leq s} \{ {\tt P}(g,q',s')\}
$$
(see \UP\  for the definition of $\Pgqs$). Let $U= V-T$. 
We show that $U$ has at most $P$ moduli maps to $M_g$.
By contradiction, suppose that for $n>P$ there exist $\psi ^1,....,\psi ^n$  different such moduli maps. Thus for every
$i$ there exists a non-isotrivial family of smooth curves $X^i\la U$ having $\psi ^i:U\la M_g$ as moduli map. There exists a
dense open subset $U'\subset U$ such that for every $u$ in $U'$ we have $\psi ^i(u) \neq \psi ^j(u) $. By the previous
discussion, there exists a smooth curve $B$ of genus at most $q$ such that the restrictions of $X^i$ to $B$ are distinct,
non-isotrivial
families of smooth curves over $B-S$, where $S:=(T\cap B)_{red}$.
Let $f_i:Y^i \la B$ be the smooth relatively minimal model of $X^i_{|B}$ over $B$. Then
$[f_i]\in \FBS$ and $[f_i]\neq [f_j]$. Since $S$ has cardinality at most $s$, 
$\FBS$ has at most $P$ elements; we reached a
contradiction.

To conclude the proof of the Lemma, fix $\psi :U\la M_g$ any such a moduli map; Lemma \UPS \  shows that there exist
at most $\Dg$ families of curves over $U$ having $\psi $ as moduli map. 
We can conclude that the number of equivalence classes of
non-isotrivial families of smooth curves 
of genus $g$ over $U$ is bounded above by $\Hgdrs = \Dg \cdot P$. $\qed$

\rk In particular we obtained, of course, that for every pair $(V,T)$, the set
$F_g(V,T)$ is finite.

\

\cl{{\bf 4. Uniformity of rational points of curves over  function fields}}

\

\ni
{\it Parshin's construction revisited.}
We  deal here with a well known construction of Parshin, the goal of which was to show that the Shafarevich
conjecture implies the Mordell conjecture.
Fix a  curve $B$, a finite subset $S\subset B$ and let $K$ be the
field of rational functions of $B$;  let $f:X\la
B$ be a non-isotrivial  fibration having smooth fiber away  from $S$. 

\

\rk
By the existence and
unicity of minimal models for (non-ruled) smooth surfaces, giving $f:X\la B$ is the same as giving the generic
fiber of it. Therefore we shall abuse notation slightly and denote by $X(K)$ the set of sections of $f$,
identifying it with the set of rational points.

\

\ni
The goal of Parshin's method is to obtain a map:
$
\alpha _X : X(K) \la \cup \FFB =:P_X
$
such that the union on the right  is over a finite set. The  construction is such that
$\a _X$  has finite fibers.
The Shafarevich conjecture  (i.e. the Parshin-Arakelov Theorem) gives that $\FBS$ is
always a finite set. Therefore $X(K)$ is finite.

We will revise Parshin's method to show that,  as $X$ and $B$
vary among curves with fixed genus $g$ and $q$, and $S$ varies among sets of fixed cardinality $s$:

\ni
(1) the fibers of $\a _X$ are uniformly bounded, 

\ni
(2) the union defining $P_X$ is over a uniformly bounded set.

\ni
By Theorem \UP  \  the sets $\FBS$ are uniformly bounded.
We shall hence conclude that  rational points  are uniformly bounded.

Let $\s \in X(K)$ be a $K$-rational point of $X$, so that  $\s$ can be viewed as a section $\s:B\la X$ of 
$f:X\la B$. Let $\Sigma =\s (B)$ be the image curve.
We are going to construct a finite covering $\rho :B'\la B$, ramified only over $S$,   a  fibration $Y'\la B'$
with a commutative diagram

$$
\matrix{Y'&\mapright {} & X\cr
\  \  \lmapdown{f'}&&\mapdown{f} \cr
B' & \mapright{\rho }& B}
$$
such that every fiber $Y'_{b'}$ is a finite covering of $X_b$  ($b={\rho (b')}$) which ramifies only over $\sigma (b)$.
The goal is to carry this out so that all numerical invariants of the new family are bounded with $q$, $g$ and $s$.
We  here provide  a brief account, see [S] or [P] for the missing details.

The section $\s$ gives a map $u:X\la \Pic ^0X/B$ defined by $u(x)=x-\s (f(x))$. Multiplication by $2$ in $\Pic ^0X/B$
yields a covering of $X$ over $B$, which, away from singular fibers, is  \'etale of degree $2^{2g}$. Let $Y$ be a
connected component of such covering, we get a new family of curves of genus $g(Y)\leq 1+2^{2g}(g-1)$; for $b\not\in S$,
each fiber
$Y_b'$  is an
\'etale covering of the corresponding fiber  $X_b$.
Let $D$ be the preimage of $\Sigma$ in $Y$. Then $D$ is \'etale of degree at most $2^{2g}$ over $B \setminus S$. 

There exists a
covering $B_1\la B$ of degree at most 
$2^{2g}(2^{2g}-1)$, ramified only over $S$  such that on the relatively minimal desingularization $Y_1$
of $B_1\times _BY$ there are two disjoint sections
$\sigma _1$ and
$\sigma _2$ with $\Sigma _i =\Im\s _i$ mapping to
$D$ via the natural morphism $Y_1\la Y$. 

One can construct a further covering $B_2 \la B_1$,
ramified only over the preimage of $S$,  such that on $Y_2 = Y_1 \times _{B_1} B_2$ the line bundle 
${\cal O}_{Y_2} (\Gamma _1 + \Gamma _2) $, given  by the
the pull-back of ${\cal }\Sigma _1+\Sigma _2$ admits a square root. 
Such a $B_2$ is obtained by first mapping $B_1$ to $\Pic ^0Y_1/B_1$ via $b \mapsto \s_1(b)-\s_2(b)$, 
and then by considering the multiplication by $2$ map $\Pic ^0Y_1/B_1\la \Pic ^0Y_1/B_1$
and defining 
$B_2= B_1 \times _{\Pic ^0Y_1/B_1}\Pic ^0Y_1/B_1$
(replacing $B_2$ by a connected component dominating $B_1$ if necessary).

This ensures 
(after a further degree-$2$ base change ramified over singular fibers)
that there exists a double covering $Y_3\la
Y_2$ having branch locus 
$\Gamma _1 + \Gamma _2$. Finally, let $Y'$ be the relatively minimal resolution of $Y_3$ over $B_2$,
let $B'=B_2$, let $\rho :B'\la B$ be the covering map and let $S'=\rho ^{-1} (S)$.
Then $f':Y'\la B'$ is the family we wanted; by
construction the fibers of $Y'$ away from $S'$ are
coverings of the fibers of $f$ of degree dividing $2^{2g+1}$, having two simple ramification points lying over
$\Sigma$. Hence $g(Y')\leq 2+2^{2g+1}(g-1)$, let $\g = 2+2^{2g+1}(g-1)$.
$\rho$ is a covering of degree bounded above by a function of only $g$ and $s$, and ramified only over $S$;
let $S'=\rho ^{-1}(S)$, then 
the pair $B',S'$ belongs to a uniformly bounded set $I$ of cardinality at most ${\tt c}(g,q,s)$
(some constant depending only on $g,\  q$ and $s$).
Moreover, the genus of $B'$ is uniformly bounded by a constant $\q $. We shall
define $\a _X(\s ) = [f':Y'\la B']$.

\proclaim Lemma \FIB. 
There exists a number $\Agqs$ such that for every  curve $B$ of
genus $q$,  for every finite subset $S\subset B$ of cardinality $s$ and for
every non-isotrivial fibration $X\la B$ of genus $g$, having smooth fiber
outside of $S$, the fibers of the map $\a _X$ have
cardinality at most
$\   \Agqs$.

\ni
{\it Proof.} Let $[f':Y'\la B']$ be an element in $\Im \a _X$ which lies in $\FFB$. As we have seen, a
$K$-rational point  $\s \in \a _X^{-1}([f'])$ determines a 
commutative diagram of objects and morphisms

$$
\matrix{Y'&\mapright {} & Y\times _B B' &\mapright {}& Y&\mapright{} &X\cr
\ \lmapdown{f'}&&\mapdown{}&&\mapdown {}&&\mapdown{f} \cr
B'&=&B' & \mapright{\rho }& B&=&B}
$$
let $\rho ':Y'\la X$ be the composition map above, then  $\sigma$ is uniquely determined as
the branch locus of
$\rho '$. Now, $\rho$ varies in a uniformly bounded set, as we have just seen. By \UDF \   so does $\rho '$.
$\qed$

\proclaim Theorem \UM .
Fix $g\geq 2 $ and $q$. There exists an integer $\Mgqs$ such that 
for every  curve $B$ of genus $q$, for every $S\subset B$ with $|S|\leq s$, and
for every  non-isotrivial curve $X$ of genus $g$, defined over $K=\C (B)$ and having
 smooth fiber outside of $S$, we have
$
\  |X(K)| \leq \Mgqs
$

 \

\pf
Maintaining the same notation as before, by the Theorem of Parshin-Arakelov,
$P_X$ is a finite set; by the previous Lemma we  get:
$
|X(K)|\leq \Agqs \cdot |P_X|.\  
$
It remains to show that $P_X$ is unformly bounded.
Recall that 
$
P_X=\cup  \FFB
$
where the union is over a uniformly bounded set. More precisely, the union is taken over 

\ni
(a) $g'\leq \g$

\ni
(b)  $(B',S')$ varies in a set  of cardinality at most ${\tt c}(g,q,s)$

\ni
(c) $q'=g(B')  \leq \q$

\ni
(d) $s'=|S'|\leq {\tt s}(g,q,s)
$

By Theorem \UP , there exists  $\Pgqs$ such that for every curve $B$ of genus $q$ and every $S$ of cardinality $s$
we have $|\FBS |\leq \Pgqs$.

Hence letting
$P^M(g,q,s):= \max  \{ {\tt P}(g',q',s')\}$ for $g'$, $q'$ and $s'$ satisfying (a) (c) and (d)  above,
we conclude
$
|P_X|\leq P^M(g,q,s) \cdot \g \cdot {\tt c}(g,q,s)
$
\  $\qed$
 
\

\ni
{\it Fields of higher transcendence degree.}
We obtain as a corollary a uniformity result about rational points of curves over function fields of any
dimension. In [Ma] it is shown that working by induction on the transcendence
degree, the Mordell conjecture for curves over function fields of any dimension 
follows from the one-dimensional case. Here we show that, by specializing to certain
curves,  from the uniform version  given by \UM \  one can obtain uniformity
statements for function fields of any dimension.  We use the same notation
introduced at the end of the previous section, to state and prove Theorem \HUP .

\proclaim Theorem \HUM .
Fix non-negative integers $g\geq 2$ and $s$ and fix a polynomial $h\in \Q[x]$. There exists a number $\Mghs$ such
that for any
$V\in \Mh$, for any closed $T\subset V$ such that $K_V^{\dim V -1}\cdot T\leq s$  and for any 
non-isotrivial curve $X$ of genus $g$, defined over $K$,  having smooth fiber away from $T$ 
we have $|X(K)|\leq \Mghs$, where $K$ is the field of functions of $V$.
 
\pf
Exactly as in the proof of \HUP , one reduces the proof to the statement below.

\proclaim Lemma \PUM.
Fix $g$, $d$, $r$ and $s$ integers; there exists a number $\Ngdrs$ such that
 for any smooth irreducible projective variety
$V\subset \P ^r$ of degree $d$, for any subvariety $T\subset V$  of degree $s$ and for any non-isotrivial
family of curves $f:X\la V$ such that for every $v\not\in T$ the fiber $X_v$ is a smooth curve
of genus $g$ , there exist at most 
$\Ngdrs$ rational sections of $f$.

\pf As in the first part of the proof of \PUP , let
$q= {\tt q }(d,r)$ be the maximum genus of a  curve of degree $d$ in $\P ^r$.
Then
$V$ is covered by smooth curves of genus at most $q$ through any of its points. 
Any such a curve, call it  $B$, intersects $T$ in at most $s$ points.

Let 
$
\Ngdrs:= \max _{q'\leq q,s'\leq s}\{ {\tt M}(g,q',s')\}
$ (see \UM )
and let us show that $f$ has at most $\Ngdrs$ sections.
By contradiction, let $n>\Ngdrs$ and suppose that there exist $\s _1,....\s _n$ distinct sections of
$f$. Observe that there exists an open, dense subset $U$ of $V$ such that $\s_i$ is
regular on it and for every $u\in U$, we have $\s _i(u)\neq \s _j(u)$ if $i\neq j$.

There exists a curve $B$ of genus at most $q$ as above, such that  $B \cap U\neq
\emptyset$ and such that the restriction 
$
f_B:X_B\la B
$
is not isotrivial. $f_B$ is a family of curves of genus $g$ having smooth fiber away from the
set $(T\cap B )_{red}$  of cardinality at most $s$.
By construction, the restrictions of $\s _1 ,...,\s _n$ to $B$ give $n$ distinct
sections of
$f_B$. This is a contradiction to \UM.
$\qed$

\

\ni
{\it Aknowledgements.} 
This research was done while I was visiting the Mathematics Departments
of the Universities of Roma 1 and
Roma 3, I am grateful to these institutions for their hospitality.
I am very thankful to Edoardo Sernesi with
whom I had  many valuable conversations on this work,
and to Cinzia Casagrande, Fabrizio Catanese, Carlo Gasbarri, Barry Mazur and Kieran O'Grady
for precious comments.

\

\ni
{\bf References}

\

\ni
[AP] A.Alzati, P.Pirola: {\it Some remarks
on the de Franchis Theorem} 
Ann.Univ.Ferrara -Sez. VII - Sc.Mat.
Vol.XXXVI, 45-52 (1990).

\ni
[A] S.Ju. Arakelov: {\it Families of algebraic curves with fixed degeneracies} Izv.
Akad. Nauk. SSSR Ser. Mat.35 (1971) no. 6, 1277-1302.

\ni
[ACGH] E.Arbarello, M.Cornalba, P.Griffiths, J.Harris: 
{\it Geometry of algebraic curves}  Grundlehren der mathematischen Wissenschaften
267 Springer.

\ni 
[B] A.Beauville: {\it Expos\'e n. 6} in  {\it  S\'eminaire sur les pinceaux des
courbes de genre au moins deux} Ast\'erisque 86 (1981).

\ni
[BV] E.Bedulev, E.Viehweg: {\it On the Shafarevich conjecture for surfaces of general type over function
fields} Preprint (1999) AG 9904124.

\ni
[CHM] L.Caporaso, J.Harris, B.Mazur: {\it Uniformity of rational points} Journal of
American Mathematics Society, Vol. 7, N. 3, January 97 
p. 1-33.

\ni
[DF] M.De Franchis: {\it Un teorema sulle involuzioni irrazionali} Rend.Circ.Mat
Palermo 36 (1913), 368.

\ni
[DM] P.Deligne, D.Mumford: 
{\it The irreducibility of the space of curves of given genus} Publ Math. IHES 36
(1969) 75-120.

\ni
[EV] H.Esnault, E.Viehweg: {\it Effecive bounds for semipositive sheaves and the height of points on curves over
complex function fields} Compositio math. 76 (1990) 69-85.

\ni
[F] G.Faltings: {\it Arakelov's Theorem for abelian varieties} Invent.Math. 73
(1983) 337-348.

\ni
[G] D.Gieseker: {\it Lectures on moduli of curves} TIFR Lecture Notes 69 (1982)
Springer.

\ni [HM] J.Harris, I.Morrison {\it Moduli of curves} Graduate texts in Math. 187
(1998) Springer

\ni
[K] E.Kani: {\it Bounds on the number of non rational subfields of a function field}
Invent. Math. 85 (1986) pp. 185-198.

\ni
[Ma] Y.Manin: {\it 
Rational points of algebraic curves over function fields} Izv. Akad. Nauk. 27 (1963),
1395-1440. 

\ni
[Mi] Y.Miyaoka: {\it Themes and variations on inequalities of Chern classes}
Sûgaku 41 (1989), no. 3, 193-207. 

\ni
[Mu] D.Mumford: {\it Stability of projective varieties} 
L'Enseignement Math\'ematique 23 (1977) p. 39-110.

\ni
[P] A.N.Parshin: {\it Algebraic curves over function fields} Izv. Akad. Nauk. SSSR
Ser. Mat.32 (1968) no. 5, 1145-1170.

\ni
[S] L.Szpiro: {\it Expos\'e n. 3} in {\it S\'eminaire sur les pinceaux des courbes
de genre au moins deux} Ast\'erisque 86 (1981).

\ni
[T] S. Tan:{\it Height inequalities of algebraic points on curves over functional
fields} J.reine angew. Math. 481 (1995) 123-135.

\ni
[V] E. Viehweg: {\it Quasi-projective moduli for polarized manifolds} Ergebnisse der
Mathematik (30) Springer.
\end